\documentclass[
    ,final            
  ]
  {aipproc}

\layoutstyle{8x11single}

\usepackage{latexsym,bm,amssymb,amsfonts,amsbsy,amsmath,graphics,graphicx,color}
\usepackage{url}

\newtheorem{theo}{Theorem}
\newtheorem{rem}{Remark}
\def\RR{\mathbb R}

\def\pmatrix{ \left( \begin{array} }
\def\endpmatrix{ \end{array} \right) }

\def\d2dxx{\frac{\partial^2}{\partial x^2}}

\def\no{\noindent}

\def\phi{\varphi}

\def\P{{\cal P}}

\def\A{{\cal A}}

\begin{document}

\title{Energy and quadratic invariants preserving integrators of Gaussian type}\thanks{Work
developed within the project ``Numerical Methods and Software for Differential
Equations''.}

\classification{65P10, 65L05.}
\keywords{Hamiltonian systems, collocation Runge-Kutta methods, symplectic
integrators, energy-preserving methods.}

\author{Luigi Brugnano}{
  address={Dipartimento di Matematica ``U.\,Dini'', Universit\`a di
Firenze, Viale Morgagni 67/A, 50134 Firenze (Italy).} }

\author{Felice Iavernaro}{
  address={Dipartimento di Matematica, Universit\`a di Bari, Via Orabona 4,
          70125 Bari (Italy).}
}

\author{Donato Trigiante}{
  address={Dipartimento di Energetica ``S.\,Stecco'', Universit\`a di
Firenze, Via Lombroso 6/17, 50134 Firenze (Italy).}
}

\begin{abstract}
Recently, a whole class of evergy-preserving integrators has been
derived for the numerical solution of Hamiltonian problems
\cite{BIT0,BIT00,BIT1}. In the mainstream of this research
\cite{BIT3}, we have defined a new family of symplectic
integrators depending on a real parameter $\alpha$ \cite{BIT4}.
For $\alpha = 0$, the corresponding method in the family becomes
the classical Gauss collocation formula of order $2s$, where $s$
denotes the number of the internal stages. For any given non-null
$\alpha$, the corresponding method remains symplectic and has
order $2s-2$: hence it may be interpreted as a $O(h^{2s-2})$
(symplectic) perturbation of the Gauss method. Under suitable
assumptions, it can be shown that the parameter $\alpha$ may be
properly tuned, at each step of the integration procedure, so as
to guarantee energy conservation in the numerical solution. The
resulting method shares the same order $2s$ as the generating
Gauss formula, and is able to preserve both energy and quadratic
invariants.
\end{abstract}

\maketitle


\section{Introduction}
When dealing with the numerical integration of canonical
Hamiltonian systems in the form
\begin{equation}\label{hamilode}
\left\{  \begin{array}{l} \dot y =  J\nabla H(y) \equiv f(y),  \\
y(t_0) = y_0 \in\RR^{2m}, \end{array} \right.
 \qquad J=\pmatrix{rr} 0 & I \\ -I & 0 \endpmatrix \in \RR^{2m \times
 2m},
\end{equation} ($I$ is the identity matrix of dimension $m$), two
main lines of investigation may be traced in the current
literature, having as objective the definition and the study of
symplectic methods and energy-conserving methods, respectively. In
fact, the symplecticity of the map and the conservation of the
energy function are the most relevant features characterizing a
Hamiltonian system.

From the very beginning of this research activity, high order
symplectic formulae were already available within the class of
Runge-Kutta methods \cite{Feng,Suris,SC}, the Gauss collocation
formulae being one noticeable example. One important implication
of symplecticity of the discrete flow, for Gauss-Legendre methods,
is the conservation of quadratic invariants. This circumstance
makes the symplecticity property of the method particularly
appealing in the numerical simulation of isolated mechanical
systems in the form \eqref{hamilode}, since it provides a precise
conservation of the total angular momentum during the time
evolution of the state vector. As a further positive consequence,
a symplectic method also conserves quadratic Hamiltonian functions
(see the monographs \cite{HLW,LR} for a thorough analysis of
symplectic methods).

Conversely, if one excludes the quadratic case, energy-conserving
methods were initially not known within classical integration
methods. The unsuccessful attempts to derive energy-preserving
Runge-Kutta methods  for polynomial Hamiltonians, culminated in the
general feeling that such methods could not even exist (see
\cite{IZ} and \cite{M2AN}). A completely new approach is represented
by {\em discrete gradient methods} which are based upon the
definition of a discrete counterpart of the gradient operator so
that energy conservation of the numerical solution is guaranteed at
each step and whatever the choice of the stepsize of integration
(see \cite{G,MQR}).

More recently, the conservation of energy has been approached by
means of the definition of the {\em discrete line integral}, in a
series of papers (such as \cite{IP1,IT3}), leading to the
definition of {\em Hamiltonian Boundary Value Methods (HBVMs)}
(see for example \cite{BIT0,BIT00,BIT1,BIT2,BIS}). They are a
class of methods able to preserve, in the discrete solution,
polynomial Hamiltonians of arbitrarily high degree (and, hence, a
{\em practical} conservation of any sufficiently differentiable
Hamiltonian. Such methods admit a Runge-Kutta formulation which
reveals their close relationship with classical collocation
formulae \cite{BIT3}. An infinity extension of HBVMs has also been
proposed in \cite{BIT2} and \cite{Ha}. These limit methods may be
interpreted as a generalization of the \textit{averaged vector
field method} defined in \cite{QMcL}.

Attempts to incorporate both symplecticity and energy conservation
into the numerical method will clash with two non-existence
results. The first \cite{GM} refers to non-integrable systems,
that is systems that do not admit  other independent first
integrals different from the Hamiltonian function itself.
According to the authors' words, it states that
\begin{quote} \em
If [the method] is symplectic, and conserved $H$ exactly, then it
is the time advance map for the exact Hamiltonian system up to a
reparametrization of time.
\end{quote}
The second negative result \cite{CFM} refers to B-series
symplectic methods applied to general (not necessarily
non-integrable) Hamiltonian systems:
\begin{quote} \em
The only symplectic method (as $B$-series) that conserves the
Hamiltonian for arbitrary $H(y)$ is the exact flow of the
differential equation.
\end{quote}
Despite these discouraging results,  in \cite{BIT4} a new class of
symplectic integrators of arbitrarily high-order has been proposed
which, under some mild assumptions (see the next section), may share
both features, in the sense  specified in the theorem below. We
prefer the use of the term ``integrator'' rather than method since,
strictly speaking, our integrator may select a different symplectic
formula from one integration step to the next, in order to enforce
the energy conservation property. In what follows, we sketch the
main ideas behind this approach. For further generalizations, as
well as for a number of numerical evidences, we refer to
\cite{BIT4}. We will begin with introducing a family of one-step
methods
\begin{equation}
\label{met_alpha} y_1(\alpha,h)=\Phi_h(y_0,\alpha)
\end{equation}
($h$ is the stepsize of integration), depending on a real
parameter $\alpha$, with the following specifics:
\begin{enumerate}
\item for any fixed choice of $\alpha \not = 0$, the corresponding method is a
symplectic Runge-Kutta method with $s$ stages and of order $2s-2$,
which exactly conserves all quadratic invariants;
\item for $\alpha=0$ one gets the Gauss collocation method (of order $2s$);
\item for any choice of $y_0$ and in a given range of the stepsize $h$,
there exists a value of the parameter, say $\alpha_0$, depending
on $y_0$ and $h$, such that $H(y_1(\alpha_0,h))=H(y_0)$ (energy
conservation).
\end{enumerate}
The parametric method (\ref{met_alpha}) realizes a
 symplectic perturbation of the Gauss method of size $O(h^{2s-2})$.
 Under suitable assumptions, as the parameter $\alpha$ ranges in a small interval centered at
 zero, the value of the numerical Hamiltonian function $H(y_1)$ will
 match $H(y(t_0+h))$ thus leading to energy conservation.
This result is formalized as follows:
\begin{theo}[Energy conservation] \label{theo1}
Under suitable assumptions, there exists a real sequence
$\{\alpha_k\}$ such that the numerical solution defined by
$y_{k+1}=\Phi_h(y_{k},\alpha_k)$, with $y_0$ defined in
(\ref{hamilode}), satisfies $H(y_k)=H(y_0)$.
\end{theo}
One important remark is in order to clarify this statement and how
it relates to the above non-existence results. Let us select the
value of the parameter $\alpha=\alpha_0$, if any, in order to
enforce the energy conservation between the two state vectors
$y_0$ and $y_1$, as indicated at item 3 above\footnote{To avoid
any misunderstanding, we emphasize that the value $\alpha_0$ is
now maintained constant, otherwise the map would fail to be
symplectic.}: the map $y \mapsto \Phi_h(y,\alpha_0)$ is symplectic
and, by definition, assures the energy conservation condition
$H(y_1)=H(y_0)$. However, it is worth noticing that it would fail
to provide a conservation of the Hamiltonian function if we
changed the initial condition $y_0$ or the stepsize $h$. For
example, in general for any $\hat y_0 \not =y_0$, we would obtain
$H(\Phi_h(\hat y_0,\alpha_0)) \not = H(y_0)$: in this case we
should change the value of the parameter $\alpha$ in order to
recover the equality condition.\footnote{More in general, the
sequence $\{\hat \alpha_k\}$ that will satisfy Theorem \ref{theo1}
starting at $\hat y_0$ will differ from the sequence
$\{\alpha_k\}$. Such sequences will be defined as the solution of
the nonlinear system (\ref{concon}), as described in the next
section.} Strictly speaking, the energy conservation property
described in Theorem~\ref{theo1} weakens the standard energy
conservation condition mentioned in the two non-existence results
stated above and hence our methods are not meant to produce a
counterexample of these statements.

\section{Definition of the methods}
Let ~$\{c_1<c_2<\dots<c_s\}$~ and ~$\{b_1,\dots, b_s\}$~ be the
abscissae and the weights of the Gauss-Legendre quadrature formula
in the interval $[0,1]$. We consider the Legendre polynomials
$P_j(\tau)$ of degree $j-1$, for $j=1,\dots,s$, shifted and
normalized in the interval $[0,1]$ so that $\int_0^1
P_i(\tau)P_j(\tau) \mathrm{d} \tau = \delta_{ij}$, for
$i,j=1,\dots,s,$ ($\delta_{ij}$ is the Kronecker symbol), and the
$s\times s$ matrix $\P = \left( P_j(c_i) \right)$. Our starting
point is the following well-known decomposition of the Butcher
array $A$ of the Gauss method of order $2s$
\cite[pp.\,77--84]{HW}:
\begin{equation}\label{A}
A= \P X_s \P^{-1},
\end{equation}
where $X_s$ is defined as
\begin{equation}\label{Xs}
X_s = \pmatrix{cccc}
\frac{1}2 & -\xi_1 &&\\
\xi_1     &0      &\ddots&\\
          &\ddots &\ddots    &-\xi_{s-1}\\
          &       &\xi_{s-1} &0\\
\endpmatrix,
\qquad\mbox{with}\qquad \xi_j=\frac{1}{2\sqrt{4j^2-1}}, \qquad
j=1,\dots,s-1.\end{equation} We now consider the matrix
$X_s(\alpha)$ obtained by perturbing \eqref{Xs} as follows:
\begin{equation}\label{Xs_alpha}
X_s(\alpha) = \pmatrix{cccc}
\frac{1}2 & -\xi_1 &&\\
\xi_1     &0      &\ddots&\\
          &\ddots &\ddots    &-(\xi_{s-1}+\alpha)\\
          &       &\xi_{s-1}+\alpha &0\\
\endpmatrix \equiv X_s + \alpha W_s,\end{equation}
where $\alpha$ is a real parameter, ~$W_s = (e_se_{s-1}^T -
e_{s-1}e_s^T)$,~ and, as usual, ~$e_j\in\RR^s$~ is the $j$th unit
vector. The family of methods (\ref{met_alpha}) we are interested
in, is formally defined by the following tableau (see
(\ref{A})--(\ref{Xs_alpha})):
\begin{equation}
\label{qgauss}
\begin{array}{c|c}\begin{array}{c} c_1\\ \vdots\\ c_s\end{array} & \A(\alpha) \\
 \hline                    &b_1\, \ldots ~ b_s
\end{array} \qquad\mbox{with}\qquad \A(\alpha)\equiv \P X_s(\alpha)\P^{-1} =  A + \alpha \P W_s
\P^{-1}.
\end{equation}
Therefore $A(0)=A$ and, moreover, the following result holds true
\cite{BIT4}.

\begin{theo} For any fixed value of $\alpha$, the Runge-Kutta method
(\ref{qgauss}) is symmetric and symplectic. For $\alpha=0$, the
usual Gauss-Legendre method of order $2s$ is recovered. For any
fixed $\alpha\ne0$, a method of order $2s-2$ is obtained.
\end{theo}

If we can choose $\alpha\equiv\alpha_0$ so that  the
energy-conservation property be satisfied, then, one obtains a
{\em (symmetric), Energy and QUadratic Invariants Preserving
(EQUIP) method}, as specified in Theorem \ref{theo1}, of Gaussian
type. Indeed, the conservation of quadratic invariants easily
follows from the structure of the matrix (\ref{Xs_alpha}) defining
the method. In conclusion, these methods will provide an exact
conservation of  all quadratic invariants, besides the Hamiltonian
function. In more details, if we denote, as usual, $Y=(Y_1^T\dots
Y_s^T)^T$ the vector of the stages, $e=(1,\dots,1)^T\in\RR^s$, and
defining the error function $g(\alpha,h)= H(y_1(\alpha,h))-
H(y_0)$, the nonlinear system, in the unknowns $Y_1,\dots,Y_s$ and
$\alpha$, that is to be solved at each step for getting energy
conservation, reads, for the given stepsize $h$,
\begin{equation}
\label{concon} \left\{ \begin{array}{l} Y = e \otimes y_0 + h
(\A(\alpha) \otimes I) F(Y), \\
g(\alpha,h) =0.
\end{array}
\right.
\end{equation}
Concerning the question about the existence of a solution of
\eqref{concon}, we make the following assumptions:
\begin{itemize}
\item[($\mathcal A_1$)] the function $g$ is analytical in a rectangle $[-\bar \alpha, \bar \alpha] \times [-\bar
h, \bar h]$ centered at the origin;
\item[($\mathcal A_2$)] let $d$ be the order of the error in the Hamiltonian function
associated with the Gauss method applied to the given Hamiltonian
system \eqref{hamilode} and the given state vector $y_0$, that is:
$$ g(0,h) = H(y_1(0,h))- H(y_0) = c_0 h^{d} + O(h^{d+1}),
\qquad c_0 \not = 0.$$ Then, we assume that for, any fixed $\alpha
\not = 0$ in a suitable neighborhood of the origin, $$g(\alpha,h)
= c(\alpha) h^{d-2} + O(h^{d-1}), \qquad c(\alpha) \not =0.$$
\end{itemize}
\begin{rem} We observe that, excluding the case where the Hamiltonian $H(y)$ is quadratic
(which would imply $g(\alpha,h) = 0$, for all $\alpha$), the error
in the numerical Hamiltonian function associated with the Gauss
method is expected to behave as $O(h^{2s+1})$. Consequently, $d\ge
2s$.\end{rem} The following result then holds true \cite{BIT4}.

\begin{theo} \label{implicit}
Under the assumptions ($\mathcal A_1$) and ($\mathcal A_2$), there
exists a function $\alpha_0=\alpha_0(h)$, defined in a
neighborhood of the origin  $(-h_0,h_0)$, such that:
\begin{itemize}
\item[(i)] $g(\alpha_0(h),h)=0$, ~for all $h\in(-h_0,h_0)$;
\qquad\qquad (ii)~ $\alpha_0(h)=\mathrm{const}\cdot h^2 + O(h^3)$.
\end{itemize}
\end{theo}

\no The next result concerns the order of convergence of the
method (\ref{concon}) (again, the proof can be found in
\cite{BIT4}).

\begin{theo}
\label{fastorder} Consider the parametric method \eqref{qgauss}
and suppose that the parameter $\alpha$ is actually a function of
the stepsize $h$, according to what stated in
Theorem~\ref{implicit}. Then, the resulting method has order $2s$.
\end{theo}

Numerical tests concerning the new EQUIP methods of Gaussian type
can be found in \cite{BIT4} and in the companion paper
\cite{BIT5}.

\bibliographystyle{aipproc}   

\begin{thebibliography}{99}

\bibitem{BIS} L.\,Brugnano, F.\,Iavernaro, T.\,Susca. Numerical comparisons
between Gauss-Legendre methods and Hamiltonian BVMs defined over Gauss points.
{\em Monografias de la Real Acedemia de Ciencias de Zaragoza} (Special Issue
devoted to the 65th birthday of Manuel Calvo). In press (2010) ({\tt
arXiv:1002.2727}).

\bibitem{BIT00}  L.\,Brugnano, F.\,Iavernaro, D.\,Trigiante.
Hamiltonian BVMs (HBVMs): a family of ``drift free'' methods for
integrating polynomial Hamiltonian problems. {\em AIP Conf. Proc.}
{\bf 1168} (2009) 715--718.

\bibitem{BIT0}  L.\,Brugnano, F.\,Iavernaro, D.\,Trigiante.
{\em The Hamiltonian BVMs (HBVMs) Homepage}, {\tt arXiv:1002.2757}.

\bibitem{BIT1} L.\,Brugnano, F.\,Iavernaro, D.\,Trigiante. Analisys
of Hamiltonian Boundary Value Methods (HBVMs): a class of energy-preserving
Runge-Kutta methods for the numerical solution of polynomial
Hamiltonian dynamical systems, {\em BIT} (2009), submitted
({\tt arXiv:0909.5659}).

\bibitem{BIT2} L.\,Brugnano, F.\,Iavernaro, D.\,Trigiante.
Hamiltonian Boundary Value Methods (Energy Preserving Discrete Line
Integral Methods). {\em Jour. of Numer. Anal. Industr. and Appl. Math.}
(2010), to appear  ({\tt arXiv:0910.3621}).

\bibitem{BIT3}  L.\,Brugnano, F.\,Iavernaro, D.\,Trigiante.
Isospectral Property of HBVMs and their connections with Runge-Kutta collocation
methods. Preprint (2010)  ({\tt arXiv:1002.4394}).

\bibitem{BIT4} L.\,Brugnano, F.\,Iavernaro, D.\,Trigiante. On the existence of
energy-preserving symplectic integrators based upon Gauss collocation formulae.
Submitted (2010) ({\tt arXiv:1005.1930}).

\bibitem{BIT5} L.\,Brugnano, F.\,Iavernaro, D.\,Trigiante.
Numerical comparisons among some methods for Hamiltonian problems.
{\em This volume}.

\bibitem{M2AN} E.\,Celledoni, R.I.\,McLachlan, D.\,McLaren, B.\,Owren,
G.R.W.\,Quispel, W.M.\, Wright. Energy preserving Runge-Kutta
methods. \textit{M2AN Math. Model. Numer. Anal.} {\bf 43} (no. 4)
(2009) 645--649.

\bibitem{CFM} P.\,Chartier, E.\,Faou, A.\,Murua. An algebraic
approach to invariant preserving integrators: the case of quadratic and
Hamiltonian invariants. {\em Numer. Math.} {\bf 103}, no.\,4  (2006)
575--590.

\bibitem{Feng} Feng Kang, Qin Meng-zhao. The symplectic methods for the
computation of Hamiltonian equations. {\em Lecture Notes in Math.} {\bf 1297}
(1987) 1--37.

\bibitem{GM} Z.\,Ge, J.E.\,Marsden. Lie-Poisson Hamilton-Jacobi
theory and Lie-Poisson integrators. {\em Phys. Lett. A}  {\bf 133} (1988)
134--139.

\bibitem{G} O.\,Gonzalez. Time integration and discrete Hamiltonian
systems. {\em J. Nonlinear Sci.} {\bf 6} (1996) 449--467.

\bibitem{Ha}  E.\,Hairer. Energy-preserving variant of collocation
methods. {\em J. Numer. Anal. Ind. Appl. Math.}  to appear.

\bibitem{HLW} E.\,Hairer, C.\,Lubich, G.\,Wanner. {\em Geometric
Numerical Integration. Structure-Preserving Algorithms for Ordinary
Differential Equations}, Second ed., Springer, Berlin, 2006.

\bibitem{HW} E.\,Hairer, C.\,Lubich, G.\,Wanner. {\em Solving
Ordinary Differential Equations II}, Second ed., Springer, Berlin,
1996.

\bibitem{IP1} F.\,Iavernaro,  B.\,Pace.  $s$-Stage Trapezoidal
Methods for the Conservation of Hamiltonian Functions of Polynomial Type.
{\em AIP Conf. Proc.} {\bf 936} (2007) 603--606.

\bibitem{IT3} F.\,Iavernaro, D.\,Trigiante. High-order symmetric
schemes for the energy conservation of polynomial Hamiltonian problems. {\em J.
Numer. Anal. Ind. Appl. Math.} {\bf 4}, no.\,1-2  (2009) 87--111.

\bibitem{IZ} A.\,Iserles, A.\,Zanna. Preserving algebraic invariants with
Runge-Kutta methods. {\em J. Comput. Appl. Math.} {\bf 125} (2000)
69--81.

\bibitem{LR} B.\,Leimkuhler, S.\,Reich. {\em Simulating Hamiltonian
Dynamics}. Cambridge University Press, Cambridge, 2004.

\bibitem{MQR} R.I.\,McLachlan, G.R.W.\,Quispel, N.\,Robidoux.
Geometric integration using discrete gradient. {\em Phil. Trans. R.
Soc. Lond. A} {\bf 357} (1999) 1021--1045.

\bibitem{QMcL} G.R.W.\,Quispel, D.I.\,McLaren. A new class of energy-preserving
numerical integration methods. {\em J. Phys. A: Math. Theor.} {\bf 41} (2008)
045206 (7pp).

\bibitem{SC} J.M.\,Sanz-Serna, M.P.\,Calvo. {\em Numerical
Hamiltonian Problems}, Chapman \& Hall, London, 1994.

\bibitem{Suris} Y.B.\,Suris. The canonicity of mappings generated by
Runge-Kutta type methods when integrating the systems $x'' =
-\partial U/\partial x$. {\em USSR Comput. Maths. Math. Phys.}
{\bf 29} (1989) 138--144.
\end{thebibliography}

\end{document}